\documentclass[a4paper]{amsart}
\usepackage[utf8]{inputenc}
\usepackage{geometry}
\usepackage{float}
\usepackage{graphicx}
\usepackage{hyperref}
\usepackage{url}
\usepackage{amsmath}
\usepackage{amssymb}
\usepackage{animate}

\usepackage[mathlines]{lineno}  

\usepackage[square,numbers]{natbib}

\title{Data assimilation: A dynamic homotopy-based coupling approach}

\author{Sebastian Reich}

\begin{document}


\maketitle

\begin{abstract}
Homotopy approaches to Bayesian inference have found widespread use especially if the Kullback-Leibler divergence between the prior and the posterior distribution is large. Here we extend one of these homotopy approach to include an underlying stochastic diffusion process. The underlying mathematical problem is closely related to the Schr\"odinger bridge problem for given marginal distributions. We demonstrate that the proposed homotopy approach provides a computationally tractable approximation to the underlying bridge problem. In particular, our implementation builds upon the widely used ensemble Kalman filter methodology and extends it to Schr\"odinger bridge problems within the context of sequential data assimilation.   
\end{abstract}


%
\section{Introduction}
%

Sequential data assimilation interlaces dynamic processes with intermittent partial state observations in order to provide reliable state estimates and their uncertainties. A wide array of numerical methods have been proposed to tackle this problem computationally. Popular methods include sequential Monte Carlo, variational inference, and  various ensemble Kalman filter formulations \cite{chopin:20,Evensenetal2022}. These methods can encounter difficulties whenever the predictive distribution is incompatible with the incoming data; in other words whenever the distance between the prior, as provided by the underlying stochastic process, and the data informed posterior distribution is large. It has long been realized that this challenge can be partially circumvented by altering the underlying stochastic process through appropriate control terms or modified proposal densities \cite{chopin:20,SR-LKNPR19,reich_data_2019,HBDD20}. Recently, the connection between devising such control terms and Schr\"odinger bridge problems \cite{CGP21} has been made explicit \cite{reich_data_2019}. However, Schr\"odinger bridge problems are notoriously difficult to solve numerically. The key contribution of this paper is to provide a computationally tracktable (sub-optimal) solution via a novel extension of established homotopy approaches \cite{daumetal2010,reich2011dynamical}. Similar to related homotopy approaches for purely Bayesian inference, the solution of certain partial differential equations (PDE) is required in order to find the desired control terms \cite{reich2011dynamical,YBM14}. In line with standard ensemble Kalman filter (EnKF) methodologies we approximate these PDEs via a constant gain approximation \cite{TdWMR17}. 
We mention that there are also alternative approaches to sequential data assimilation or inference which utilize ideas from optimal transportation; see for example \cite{SR-R13,corenflos2021differentiable,spantini2019coupling,BHDJ19}.

The paper is structured as follows. The mathematical formulation of the data assimilation problem, as considered in this paper,  is laid out in Section \ref{sec:Problem}. The standard optimal control and Schr\"odinger bridge approach to data assimilation is briefly summarized in Section \ref{sec:SB}, and the novel control formulation based on an homotopy formulation is introduced in Section \ref{sec:Homotopy}. A practical implementation based on the EnKF methodology is proposed in Section \ref{sec:Implementation}. A series of increasingly complex data assimilation problems is considered in Section \ref{sec:Examples} in order to demonstrate the feasibility of the proposed methodologies. The paper concludes with some conclusions in Section \ref{sec:Conclusions}.
Detailed mathematical derivations can be found in Appendix \ref{sec:forcing_derivation} and \ref{sec:linear_moment_matching}, respectively. 

%
\section{Problem Formulation and Background}
\label{sec:Problem}
%

We consider drift diffusion processes given by a
stochastic differential equations (SDE)
\begin{equation}
    \label{base_sde}
    {\rm d}X_t = f(X_t){\rm d}t + \sqrt{2\sigma}{\rm d}W_t,
\end{equation}
where $X_t : \Omega \to \mathbb{R}^{d_x}$, $f:\mathbb{R}^{d_x} \to \mathbb{R}^{d_x}$, $\sigma \in \mathbb{R}_{\geq 0}$, and $W_t:\Omega \to \mathbb{R}^{d_x}$ denotes
$d_x$-dimensional standard Brownian motion \cite{oksendal2013stochastic,pavliotis_stochastic_2014}.

Assuming the law of $X_t$ is absolutely continuous w.r.t Lebesgue measure with density $\pi_t$, this leads to the Fokker-Planck equation \cite{pavliotis_stochastic_2014}
\begin{equation}
    \label{fp}
    \partial_t \pi_t = -\nabla  \cdot \left( \pi_t \left( f - \sigma \nabla \log \pi_t \right) \right) .
\end{equation}
The SDE (\ref{base_sde}) can be replaced by the mean field ODE
\begin{equation}
    \label{base_ode}
    \frac{{\rm d}}{{\rm d}t} \tilde X_t = f(\tilde X_t) -
    \sigma \nabla \log \tilde \pi_t
\end{equation}
where $\tilde \pi_t$ denotes the law of $\tilde X_t$. Provided
$\tilde \pi_0 = \pi_0$, it holds that $\tilde \pi_t = \pi_t$ for all
$t> 0$. Note that the evolution of the random variable $\tilde X_t$ is entirely deterministic subject to random initial conditions $\tilde X_0 \sim \pi_0$.

At time $t = T>0$, we have observations of the system according to
\begin{equation}
    \nonumber
    y_T = h(x_T^\dagger) + \nu
\end{equation}
from which we wish to infer the unknown state $x_T^\dagger \in \mathbb{R}^{d_x}$. Here $h:\mathbb{R}^{d_x} \to \mathbb{R}^{d_y}$ denotes the forward map and $\nu \sim \mathcal{N}(0,R)$ is $d_y$-dimensional Gaussian noise with covariance matrix $R \in \mathbb{R}^{d_y \times d_y}$.

Let $L:\mathbb{R}^{d_x} \to \mathbb{R}$ denote the corresponding negative log-likelihood function. Since $\nu$ is Gaussian it is given by
\begin{equation}
    \label{likelihood}
    L(x) = \frac{1}{2}(h(x) - y_T)^\top R^{-1} (h(x) - y_T)
    .
\end{equation}
The observations are combined with the predictive density $\pi_T$ at time $t=T$ according to Bayes' theorem,
\begin{equation} \label{Bayes}
    \pi^{\rm a}_T = \frac{e^{-L} \pi_T}{\int e^{-L(x)} \pi_T(x) {\rm d}x} .
\end{equation}
The process of transforming the random variable $X_T \sim \pi_T$ into a random variable $X_T^a \sim \pi_T^{\rm a}$ is called data assimilation
in the context of dynamical systems and stochastic processes \cite{law2015data,reich2015probabilistic,Evensenetal2022}.

Since performing data assimilation can be difficult if the relative Kullback-Leibler divergence
$$
{\rm KL}(\pi_T|\pi_T^{\rm a}) = \int_{\mathbb{R}^{d_x}} \pi_T(x) 
(\log \pi_T(x)-\log \pi_T^{\rm a}(x)) {\rm d}x,
$$
also called the relative entropy \cite{pavliotis_stochastic_2014}, between the prior $\pi_T$ and posterior $\pi_T^a$ is large and/or if the involved distributions are strongly non-Gaussian \cite{snyder2008obstacles,agapiou2017importance}, we propose to construct a new SDE with state process $X_t^{\rm h}$ such that $X_0^{\rm h} \sim \pi_0$ and $X_T^{\rm h} \sim \pi_T^a$. In other words, we are looking for a stochastic process (bridge) with initial density $\pi_0$ and final density $\pi^{\rm a}_T$. The problem of finding the optimal process (in the sense of minimal Kullback-Leibler divergence) is known as the Schrödinger bridge problem \cite{CGP21}. 

%
\section{Schr\"odinger Bridge Approach}
\label{sec:SB}
%

The Bayesian adjustment (\ref{Bayes}) at final time $t=T$ leads in fact to an adjustment over the whole solution space of the underlying diffusion process described by (\ref{base_sde}). Let us denote the so called smoothing distribution by $\pi_t^{\rm a}$, $t \in [0,T]$ \cite{sr:sarkka,chopin:20}. It is well established that these marginal distributions can be generated from a controlled SDE
\begin{equation} \label{smoothing SDE}
{\rm d}X_t^{\rm a} = f(X^{\rm a}_t){\rm d}t 
    + g_t^{\rm a}(X^{\rm a}_t) {\rm d}t + \sqrt{2\sigma} {\rm d}W_t
\end{equation}
for appropriate control $g_t^{\rm a}:\mathbb{R}^{d_x} \to \mathbb{R}^{d_x}$ such that $X_0^{\rm a}\sim \pi_0^{\rm a}$ implies
$X_t^{\rm a}\sim \pi_t^{\rm a}$ for all $t>0$. It is also well known that finding a suitable $g_t^{\rm a}$ can be formulated as an optimal control problem which in turn is closely related to the backward Kolmogorov equation \cite{pavliotis_stochastic_2014,reich_data_2019}. Formulations related to (\ref{smoothing SDE}) have also been used in the context of sequential Monte Carlo methods \cite{chopin:20}.

As proposed in \cite{reich_data_2019}, an alternative perspective on sequential data assimilation is provided by Schr\"odinger bridges.
Given two marginal distributions $q_0$ and $q_T$ and stochastic process $X_t$ (referred to as the \emph{reference process}), a Schrödinger bridge is another stochastic process $\hat{X}_t$ such that $\hat{X}_0 \sim q_0$, $\hat{X}_T \sim q_T$ and the Kullback-Leibler divergence between the processes $\{\hat{X}_t\}_{t\in [0,T]}$ and $\{X_t\}_{t\in [0,T]}$ is minimal. In the language we have been using so far that means the marginals are the initial and posterior densities, i.e. $q_0 = \pi_0$ and $q_T = \pi^{\rm a}_T$, and the reference process is the solution to \eqref{base_sde}. The solution to the associated Schr\"odinger bridge problem is again of the form (\ref{smoothing SDE}) with modified control term denoted by $g_t^{\rm SB}(x)$. 

A Schr\"odinger bridge is thus the optimal coupling as measured by the Kullback-Leibler divergence  to the underlying reference process. Unfortunately the required control term $g_t^{\rm SB}$ seems rather difficult to compute in practice. In addition to the computational complexity of solving nonlinear Schr\"odinger bridge problems, the target distribution $\pi_T^a$ is implicitly defined in the setting of data assimilation. The next section offers a solution to both of these issues. We point to \cite{SR-LKNPR19} for a discussion of alternative approaches which introduce appropriate control terms into data assimilation procedures. 

%
\section{Homotopy Induced Dynamic Coupling}
\label{sec:Homotopy}
%

Since Schr\"odinger bridges are computationally challenging, we ask whether a less optimal but cheaper approach might also be feasible. Indeed, in the context of data assimilation a non-optimal coupling can be found via a homotopy between the initial and target distribution as follows. Let
\begin{equation}
    \nonumber
    \pi^{\rm h}_t(x) = Z^{-1}_t e^{-\frac{t}{T}L(x)} \pi_t(x)
\end{equation}
denote the homotopy in question, with $Z_t = \int e^{-\frac{t}{T}L(x)} \pi_t(x) {\rm d}x$ the time dependent normalization constant. It clearly holds that $\pi^{\rm h}_0 = \pi_0$ and $\pi^{\rm h}_T = \pi_T^{\rm a}$.
Note that the scaling $t \mapsto e^{-\frac{t}{T}L}$ was chosen for its simplicity and follows previous work on Bayesian inference problems \cite{daumetal2010,reich2011dynamical}. Finding better homotopies or systematic ways of constructing one could be an interesting direction for future research.

We can then reason backwards from the Fokker-Planck equation of $\pi^h_t$ to conclude that if it is the density of a random variable $X^h_t$ then that random variable must satisfy the modified SDE:
\begin{equation} \label{control_sde1}
    {\rm d}X_t^{\rm h} = f(X^{\rm h}_t){\rm d}t 
    - \frac{\sigma t}{T} \nabla L(X^{\rm h}_t){\rm d}t 
    + g_t(X^{\rm h}_t) {\rm d}t + \sqrt{2\sigma} {\rm d}W_t ,
\end{equation}
where $g_t$ is a solution to the PDE
\begin{equation}
    \label{g_t}
    \nabla \cdot (\pi^{\rm h}_t g_t) =
    \frac{1}{T} \pi^h_t \left( L + t \nabla L \cdot
        \left( f - \frac{\sigma t}{T} \nabla L - \sigma \nabla \log \pi^{\rm h}_t
    \right)
    \right)
    + \pi^{\rm h}_t \frac{\dot{Z}_t}{Z_t} .
\end{equation}
The derivations of \eqref{g_t} can be found in Appendix \ref{sec:forcing_derivation}. 

Note that (\ref{control_sde1}) constitutes a mean field model since $g_t$ depends on the distribution $\pi_t^{\rm h}$ of $X_t^{\rm h}$. We also wish to point out that
\begin{equation}
    \hat g_t^{\rm SB}(x) := - \frac{\sigma t}{T} \nabla L(x)
    + g_t(x)
\end{equation}
provides a (non-optimal) solution to the associated Schr\"odinger bridge problem. 

Since \eqref{g_t} is linear in $g_t$ we can decompose (\ref{g_t}) into a set of simpler equations $\nabla \cdot \left( \pi^{\rm h}_t g^i_t \right) = \pi^{\rm h}_t (k^i - \mathbb{E} k^i)$ such that the $k^i$ add up to the right hand side of \eqref{g_t}. In order to maintain $\int_{}^{} \nabla \cdot \left( \pi^{\rm h}_t g_t \right) = 0$ for the individual $g^i_t$ we can make use of the fact that the terms in \eqref{g_t} are of the form $\pi^{\rm h}_t \left( k - \mathbb{E} k \right)$. Separating the terms, we obtain the following equations, the sum of whose solutions solves \eqref{g_t}:
\begin{subequations}
\begin{align}
    \label{sep_g1}
    \nabla \cdot \left( \pi^{\rm h}_t g_t^1 \right) &= \pi^{\rm h}_t \left(  \frac{L}{T} - \mathbb{E}\frac{L}{T} \right) \\
    \label{sep_g2}
    \nabla \cdot \left( \pi^{\rm h}_t g_t^2 \right) &= \pi^{\rm h}_t \left( \frac{t}{T} \nabla L \cdot f - \mathbb{E} \frac{t}{T} \nabla L \cdot f  \right) \\
    \label{sep_g3}
    \nabla \cdot \left( \pi^{\rm h}_t g_t^3 \right) &= - \pi^{\rm h}_t \left( \frac{\sigma t}{T} \nabla L \cdot \nabla \log \pi^{\rm h}_t - \mathbb{E} \frac{\sigma t}{T} \nabla L \cdot \nabla \log \pi^{\rm h}_t \right) \\
    \label{sep_g4}
    \nabla \cdot \left( \pi^{\rm h}_t g_t^4 \right) &= - \pi^{\rm h}_t \left( \frac{\sigma t^2}{T^2} \nabla L \cdot \nabla L - \mathbb{E} \frac{\sigma t^2}{T^2} \nabla L \cdot \nabla L \right)
    .
\end{align}
\end{subequations}

Note that 
\begin{equation} \nonumber
\pi^{\rm h}_t \nabla L \cdot \nabla \log \pi^{\rm h}_t  = \nabla \left( \pi^{\rm h}_t \nabla L \right) - \pi^{\rm h}_t \Delta L ,
\end{equation}
which can be used to avoid the computation of $\nabla \log \pi_t^{\rm h}$ (with $\Delta = \nabla \cdot \nabla$ the Laplacian operator). Thus the controlled SDE (\ref{control_sde1}) can be replaced by
\begin{equation} \label{control_sde2}
{\rm d}X_t^{\rm h} = f(X^{\rm h}_t){\rm d}t 
    - \frac{2\sigma t}{T} \nabla L(X^{\rm h}_t){\rm d}t 
    + \hat g_t(X^{\rm h}_t) {\rm d}t + \sqrt{2\sigma} {\rm d}W_t ,
\end{equation}
where $\hat g_t$ is a solution to
\begin{equation}
    \label{hatg_t}
    \nabla \cdot (\pi^{\rm h}_t \hat g_t) =
    \frac{1}{T} \pi^{\rm h}_t \left( L + t \nabla L \cdot
        \left( f - \frac{\sigma t}{T} \nabla L 
    \right) + \sigma t \Delta L
    \right)
    + \pi^{\rm h}_t \frac{\dot{Z}_t}{Z_t} .
\end{equation}
Furthermore, if $\Delta L$ is a constant (as a function of $x$) or small in comparison to the other contributions in (\ref{hatg_t}), then
(\ref{hatg_t}) simplifies further. In particular, this is the case if the forward map is linear, that is, $h(x) = Hx$.

Building upon the mean field ODE (\ref{base_ode}), one obtains the equivalent controlled mean field ODE system
\begin{equation} \nonumber
    \frac{{\rm d}}{{\rm d}t}\tilde X_t^{\rm h} = f(\tilde X^{\rm h}_t)
    - \frac{\sigma t}{T} \nabla L(\tilde X^{\rm h}_t) - \sigma\nabla \log \tilde \pi^{\rm h}_t(\tilde X_t^{\rm h})
    + g_t(\tilde X^{\rm h}_t)  ,
\end{equation}
with $g_t$ defined as before. This mean field formulation again requires knowledge (or approximation) of $\nabla \log \tilde \pi_t^{\rm h}$. A Gaussian approximation might be sufficient in certain circumstances giving rise to
$$
\nabla \log \tilde \pi_t^{\rm h}(x) \approx -(\tilde \Sigma_t^{
\rm h})^{-1}(x-\tilde \mu_t^{\rm h})
$$
where $\tilde \mu_t^{\rm h}$ denotes the mean of $\tilde X_t^{\rm h}$ and $\tilde \Sigma_t^{\rm h}$ its covariance matrix.

%
\section{Numerical implementation}
\label{sec:Implementation}
%

No analytic solution to \eqref{g_t} is known and we thus have to resort to approximations. We note that a similar PDE arises in the computation of the gain in the feedback particle filter \cite{yang_feedback_2013} and one could use the diffusion map based approximation \cite{taghvaei_diffusion_2019} for the problem at hand. This method also transforms the PDE into a Poisson equation which it then translated into an equivalent integral equation, the semi-group form of the Poisson equation. As the name suggests the integral equation makes use of the generator of a semi-group which can be approximated by diffusion maps. Here we instead propose to follow the constant gain approximation first introduced in the EnKF methodology \cite{TdWMR17}.

%
\subsection{Ensemble Kalman mean field approximation}
%

Let us assume that $\Delta L \approx const$ in (\ref{hatg_t}). Then we only need to deal with the modified negative log likelihood function
\begin{align} \nonumber
\tilde L(x) &= \frac{1}{T} \left( L(x) + t\nabla L(x) \cdot \left(f(x)-\frac{\sigma t}{T} \nabla L(x)\right)\right) \\ \nonumber
&\approx \frac{1}{T} \left( L(x) + \frac{t}{\Delta t}\left\{
L(x) - L\left( x-\Delta t f(x) + \Delta t\frac{t\sigma}{T} \nabla L(x)\right)\right\} \right)
\end{align}
with $\Delta t$ being the time-step also used later for time-stepping the evolution equations (\ref{control_sde1}) or
(\ref{control_sde2}), respectively. Since $L$ is given by (\ref{likelihood}), we define the modified forward map
$$
\tilde h(x) = h\left(x - \Delta t f(x) + \Delta t \frac{\sigma t}{T}
\nabla L(x)\right)
$$
and thus
\begin{equation} \label{mod_likelihood}
\tilde L(x) \approx \frac{t+\Delta t}{2\Delta t T}
(h(x)-y_T)^\top R^{-1}(h(x)-y_T) - \frac{t}{2\Delta t T}
(\tilde h(x)-y_T)^\top R^{-1}(\tilde h(x) - y_T).
\end{equation}
Following the standard EnKF methodology, this suggests to approximate
the drift function $\hat g_t$ in (\ref{control_sde2}) as follows:
\begin{equation} \label{control_constant_gain}
    \hat g_t^{\rm KF}(x)= -\frac{t + \Delta t}{\Delta t T}\Sigma_t^{xh}R^{-1}\left( \frac{1}{2} \left(
    h(x)+ \pi_t^{\rm h}[h]\right) - y_T\right)
    + \frac{t}{\Delta t T}\Sigma_t^{x\tilde h}R^{-1}
    \left(\frac{1}{2} \left(\tilde h(x) +\pi_t^{\rm h}[\tilde h] \right) -y_T\right).
\end{equation}
Here we have introduced the notation $\pi_t^{\rm h}[l]$ to denote the expectation value $\mathbb{E} l$ of a function $l(x)$ under the PDF $\pi_t^{\rm h}$. Furthermore, $\Sigma_t^{xh}$ denotes the correlation matrix between $x$ and $h(x)$ under the PDF $\pi_t^{\rm h}$ etc. The derivation of (\ref{control_constant_gain}) can be found in Appendix \ref{sec:linear_moment_matching}.

%
\subsection{Particle approximation and time-stepping}
%

The controlled mean field equations (\ref{control_sde2}) can be implemented numerically by the standard Monte Carlo {\it ansatz}, that is, $M$ particles $X_t^{(i)}$ are propagated according to
\begin{equation} \label{Monte Carlo}
{\rm d}X_t^{(i)} = f(X^{(i)}_t){\rm d}t 
    - \frac{2\sigma t}{T} \nabla L(X^{(i)}_t){\rm d}t 
    + \hat g_t^{\rm KF}(X^{(i)}_t) {\rm d}t + \sqrt{2\sigma} {\rm d}W_t^{(i)}
\end{equation}
for $i=1,\ldots,M$. The required expectation values in $\hat g_t^{\rm KF}$ are evaluated with respect to the empirical measure
$$
\hat \pi_t^{\rm h}(x) = \frac{1}{2} \sum_{i=1}^M \delta (x-X_t^{(i)}).
$$
The interacting particle system can be time-stepped using an appropriate adaptation of (\ref{time stepping}) from Appendix \ref{sec:linear_moment_matching}. The computation of gradients can be avoided by applying the statistical linearisation (\ref{Stein}).

%
\section{Examples}
\label{sec:Examples}
%

We now discuss a sequence of increasingly complex examples. The purpose is both to illuminate certain aspects of the proposed control terms as well as to indicate the computational advantages of the proposed methodology. All examples will be based on linear forward maps $h(x) = Hx$ and, therefore $\Delta L$ is constant and can be ignored.

%
\subsection{Pure diffusion processes}
%

We set the drift $f$ to zero in (\ref{base_sde}) and also assuming Gaussian initial conditions. Then the homotopy approach gives rise to the controlled SDE
\begin{subequations} \label{pure_diffusion}
\begin{align}
    {\rm d}X_t^{\rm h} &= \sqrt{2\sigma}{\rm d}W_t -
    \frac{2\sigma t}{T} H^\top R^{-1} (H X_t^{\rm h}-y_T) {\rm d}t \\
    & \qquad 
    -\,\,\Sigma_t^{\rm h}H^\top  \left\{ 
    \frac{1}{T} R^{-1} - \frac{2\sigma t^2}{T^2} R^{-1} H H^\top R^{-1} \right\}
       \left(\frac{1}{2} \left(HX_t^{\rm h}  + H\mu_t^{\rm h}
     \right)-y_T  \right){\rm d}t .
\end{align}
\end{subequations}
We note that $X_t^{\rm h} \sim \pi_t^{\rm h}$ will remain Gaussian for all times and we denote the mean by $\mu_t^{\rm h}$ and the covariance matrix by $\Sigma_t^{\rm h}$. Hence, it holds that $\Sigma_t^{xx} = \Sigma_t^{\rm h}$ and $\pi_t^{\rm h}[x] = \mu_t^{\rm h}$. 

Please note that the additional drift term in (\ref{pure_diffusion}a) is pulling $X_t^{\rm h}$ towards the observation $y_T$ regardless of the value of $\Sigma_t^{\rm h}$. It should also be noted that the drift term in (\ref{pure_diffusion}b) can be both attractive or repulsive with regard to the observation $y_T$ depending on the eigenvalues of
\begin{equation} \nonumber
\Omega(t) = \frac{1}{T} R^{-1} - \frac{2\sigma t^2}{T^2} R^{-1} H H^\top R^{-1}.
\end{equation}
The strength of this drift term is moderated by the covariance matrix 
$\Sigma_t^{\rm h}$.

We consider a one-dimensional problem with $R = 0.01$, $\sigma = 1$, $H = 1$, $y_T = 1$ and $T=1$. The initial conditions are Gaussian with mean $\mu_0 = 0$ and variance $\Sigma_0 = 1$. It follows that $\pi_1$ is Gaussian with mean $\mu_1 = 0$ and variance $\Sigma_1 = 2$ and the resulting Gaussian posterior $\pi_1^a$ has mean and variance given by
\begin{equation} \label{eq:post_diffusion}
    \mu_1^a = K y_T \approx 0.9524, \qquad \mu_1^a = 2 - 2K \approx 0.0952
\end{equation}
with Kalman gain $K = 2/(2+0.01) \approx 0.9524$. 

\begin{figure}[ht]
    \centering
    \includegraphics[width=0.48\textwidth]{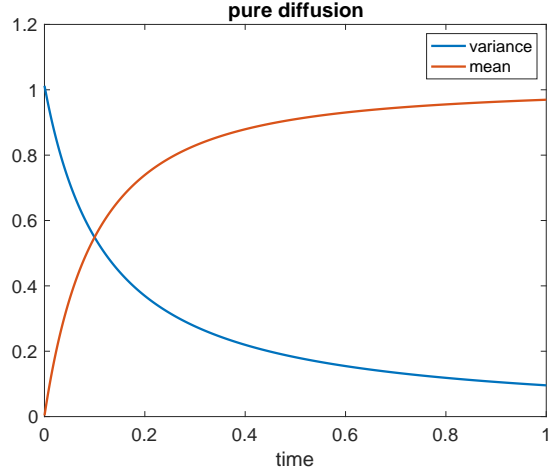}
    \caption{Time evolution of the mean $\mu_t^h$ and the variance $\Sigma_t^h$ under the mean field equations (\ref{pure_diffusion}). Their values at final agree with the posterior values provided by (\ref{eq:post_diffusion}). 
    }
    \label{fig:diffusion}
\end{figure}

In Figure \ref{fig:diffusion} one can find the time evolution of the mean and the variance under the mean field equations (\ref{pure_diffusion}). The early impact of the data driven control term on the dynamics is perhaps surprising and quite opposite to the standard sequential approach to data assimilation where one first propagates to final time and only then adjusts according to the available data.
It is also worth noticing that the sign of the corresponding 
$\Omega(t)$ changes sign at $t_c = \sqrt{2}/20$ implying that the drift term in (\ref{pure_diffusion}) has a destabilizing effect on the dynamics for $t>t_c$.

%
\subsection{Purely deterministic processes}

We now set $\sigma =0$ in (\ref{base_sde}). 
We obtain from (\ref{control_constant_gain}) the mean field ODE system
\begin{subequations}
\begin{align}
    \frac{\rm d}{{\rm d}t} X_t^{\rm h} &= f(X_t^{\rm h}) -
    \frac{t+\Delta t}{\Delta t T} \Sigma_t^{\rm h} H^\top 
    R^{-1} \left( \frac{1}{2} H \left(X_t^{\rm h} + \mu_t^{\rm h}\right) - y_T \right) \\
    & \qquad +\,\,\frac{t}{\Delta t T} \Sigma_t^{x\tilde h} R^{-1}
    \left( \frac{1}{2} \left( \tilde h(X_t^{\rm h}) + \pi_t^{\rm h}[\tilde h]\right) - y_T \right)
\end{align}
\end{subequations}
with 
\begin{equation} \nonumber
\tilde h(x) = Hx - \Delta t H f(x).
\end{equation}
These equations can be expanded giving rise to
\begin{subequations} \label{pure_drift}
\begin{align}
    \frac{\rm d}{{\rm d}t} X_t^{\rm h} &= f(X_t^{\rm h}) -
    \frac{1}{T} \left\{ \Sigma_t^{\rm h} + t \Sigma_t^{xf}\right\} H^\top 
    R^{-1} \left( \frac{1}{2} H \left(X_t^{\rm h} + \mu_t^{\rm h}\right) - y_T \right) \\
    & \qquad -\,\,\frac{t}{2T} \Sigma_t^{\rm h} H^\top R^{-1}
    H \left( f(X_t^{\rm h}) + \pi_t^{\rm h}[f])\right) 
\end{align}
\end{subequations}
upon ignoring terms of order $\mathcal{O}(\Delta t)$. Unless the drift function $f$ is linear, these mean field equations provide only an approximation to the controlled mean field equations (\ref{control_sde2}).

%
\subsection{Linear Gaussian case} \label{example:Gauss}
%

It is instructive to investigate the linear case
\begin{equation} \label{drift_linear}
f(x) = Fx + b
\end{equation}
in more detail where again everything remains Gaussian provided 
$X_0^h$ is Gaussian distributed, that is $\pi_0(x) = \mathcal{N}(x; \mu_0, \Sigma_0)$. Under these conditions the densities $\pi_t$ and $\pi^{\rm h}_t$ will also be Gaussian, we write $\pi^{\rm h}_t(x) = \mathcal{N}(x; \mu^{\rm h}_t, \Sigma^{\rm h}_t)$. The associated mean field equations follow from Appendix \ref{sec:linear_moment_matching} and are given by
\begin{subequations} \label{linear_homotopy}
\begin{align}
    \frac{{\rm d}}{{\rm d}t} X_t^{\rm h} &= FX_t^{\rm h} + b + \sigma (\Sigma_t^{\rm h})^{-1}
    (X_t^{\rm h}-\mu_t^{\rm h}) - \frac{2\sigma t}{T} H^\top R^{-1} (H X_t^{\rm h}-y_T)  \\
    & \quad -\,\,
    \left\{ \frac{1}{T} \Sigma_t^{\rm h} + \frac{t}{T} 
    \Sigma_t^{\rm h} F^\top - \frac{2\sigma t^2}{T^2} 
    \Sigma_t^{\rm h} H^\top R^{-1} H \right\} H^\top 
    R^{-1} \left( \frac{1}{2} H \left(X_t^{\rm h} + \mu_t^{\rm h}\right) - y_T \right)\\ & \quad -\,\,\frac{t}{T} \Sigma_t^h H^\top R^{-1}
    H \left( \frac{1}{2} F\left( X_t^{\rm h} + \mu_t^{\rm h} \right) + b\right).
\end{align}
\end{subequations}

A qualitative discussion can be performed in the scalar case, that is  $d_x = 1$, $H = 1$, $\sigma = 1$, $b=0$, $T=1$ and $F = \lambda$. One finds that the control terms involving $F$ stabilize the dynamics whenever $\lambda >0$. This observation is in line with the fact that the data is crucial only if the dynamics in $X_t$ is unstable, that is, $\lambda > 0$. 

%
%

We consider a two dimensional diffusion process with state variable $x = (x_1,x_2)^\top$ and linear drift term (\ref{drift_linear}) given by
\begin{equation} \nonumber
F =
\big( \begin{smallmatrix}
    -2 & 1 \\
    1 & -2
\end{smallmatrix} \big)
,
\end{equation}
$b = 0$, and diffusion constant $\sigma = 0.1I$. The forward operator is $H=
\begin{pmatrix}
    1 & 0
\end{pmatrix}
$
and the variance of the noise $R=0.01$.
The initial distribution was $\pi_0 = \mathcal{N}( (1, 3), 0.02I )$.
The observed value at time $T=1$ is set to $y_T = 2.5$. 

The posterior mean takes values $\mu_1^{\rm a} \approx 2.25$ and $\mu_2^{\rm a} \approx 1.50$, while the posterior covariance matrix
becomes
$$
\Sigma^{\rm a} \approx \left( \begin{array}{cc} 0.0086 & 0.0039\\
0.0039 & 0.0503 \end{array} \right).
$$

\begin{figure}[ht]
    \centering
    \includegraphics[width=0.48\textwidth]{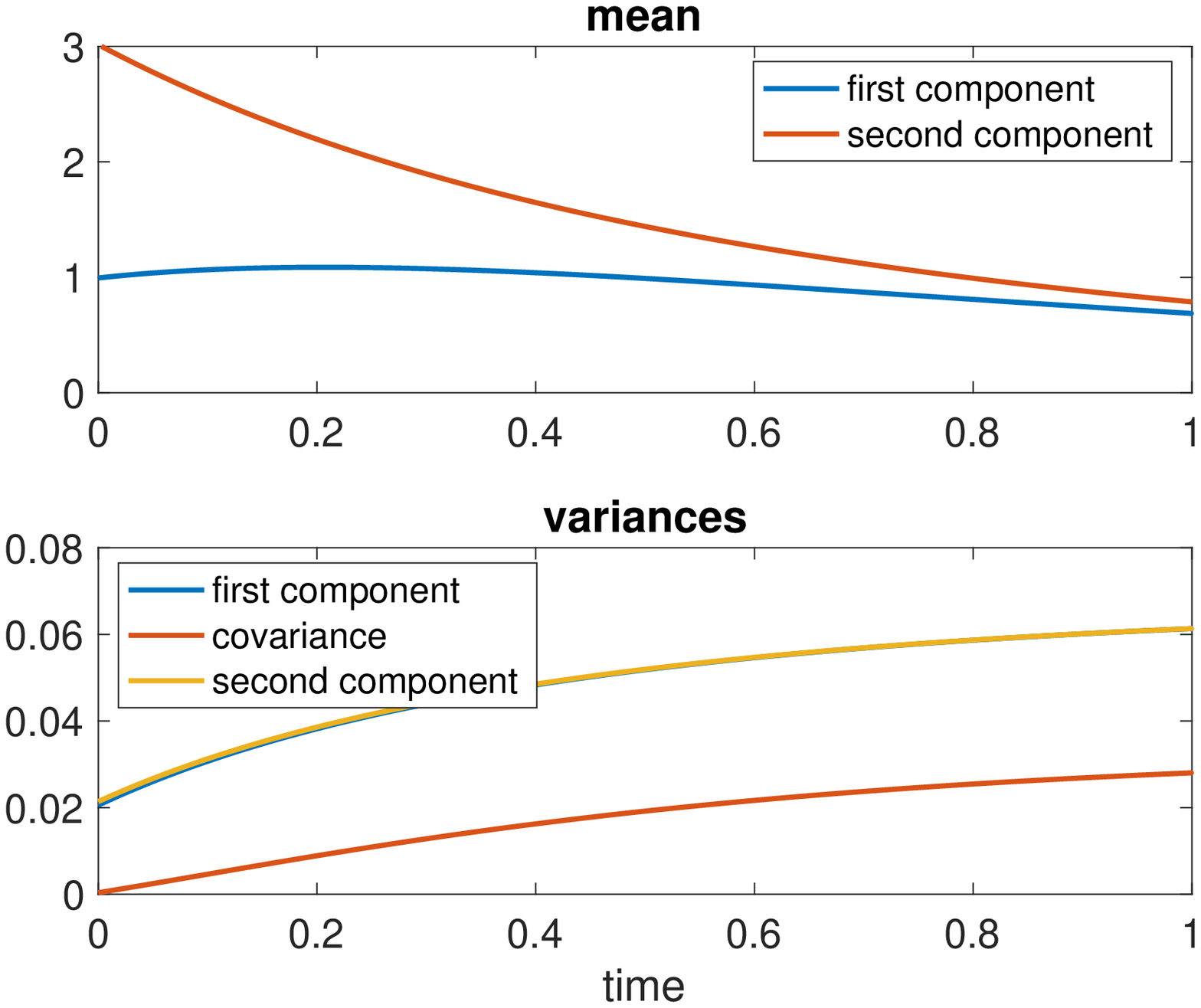}
    $\quad$
    \includegraphics[width=0.48\textwidth]{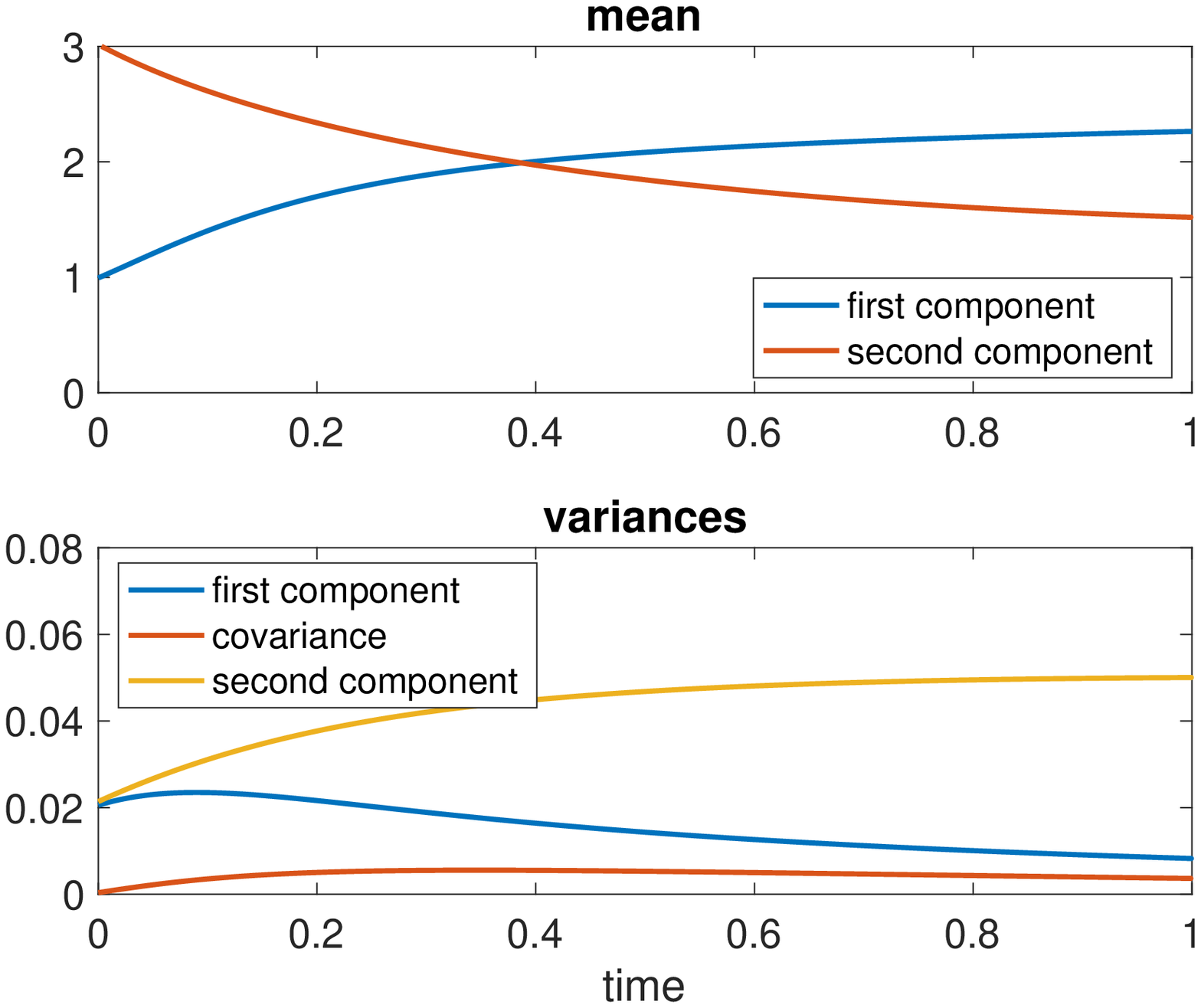}
    \caption{Left panel: Time evolution of the mean in $x_1$ and $x_2$, the two associated variances and the covariance between $x_1$ and $x_2$ under the linear diffusion process. Right panel: 
    Time evolution of the the same quantities under the controlled diffusion process.
    }
    \label{fig:twoDim}
\end{figure}

Numerical results can be found in Figure \ref{fig:twoDim}. The impact of the control term on the linear diffusion process can clearly be seen and is most prominent on the observed $x_1$ component of the process. The final values of the controlled process agree well with their posterior
counterparts.

%
%

%

%
\subsection{Nonlinear diffusion example}
%

We consider a two-dimensional problem and denote the state variable by
$x = (x_1,x_2)^\top$. The drift term is given by
$$
f(x) = -\nabla V(x), \qquad V(x) = 
\frac{\lambda_1}{2} \left(x_2 - 2 + \beta x_1^2\right)^2 +
\frac{\lambda_2}{2}\left(\frac{x_1^4}{2}-x_1^2\right) 
$$
with parameters $\lambda_1 = 2000$, $\lambda_2 = 5$, and $\beta = 1/5$. The diffusion constant is set to $\sigma = 1$. The choice of the potential $V(x)$
has two effects: (i) there is a relative high barrier for particles to pass from positive to negative $x_1$-values and vice versa; (ii) the dynamics stay close to the parabola $x_2 = 2- \beta x_1^2$.

The initial distribution is obtained by sampling $x_1$ from a Gaussian with mean
$1.5$ and variance $0.0625$. The $x_2$ component is obtained from the relation
$$
x_2 = 2 - \beta x_1^2.
$$

We observe the first component $x_1$ of the state vector at time $T=1$ with measurement error variance $R = 0.01$. The observed value is set to $y_T = -1.5$.
Due to the tiny observation error the posterior is centered sharply about the observed value. Furthermore, recall that the dynamics is essentially slaved to the parabola $x_2 = 2 - \beta x_1^2$ which makes the inference problem strongly nonlinear.

All particle simulations are run with an ensemble size of $M=1000$. Essentially identical results are obtained for $M=100$. Smaller ensemble sizes lead to numerical instabilities. 

\begin{figure}[ht]
    \centering
    \includegraphics[width=0.48\textwidth]{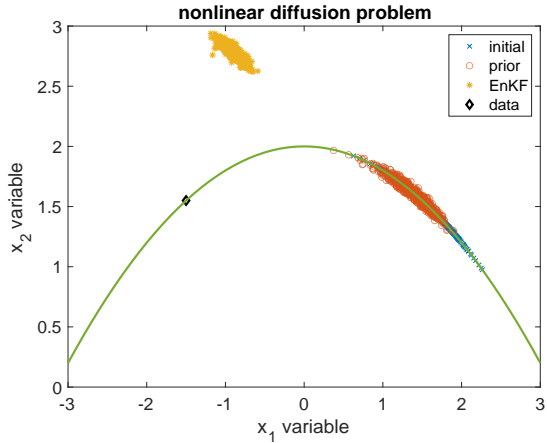}
    \caption{
        Initial (blue) and final particle positions (red) under the given evolution process together with the posterior approximation provided by the EnKF (yellow). The observed value is also displayed.
    }
    \label{fig:nonlinear1}
\end{figure}

In Figure \ref{fig:nonlinear1}, one can find the particle distribution at time $t=1$ which constitutes the prior distribution for the associated Bayesian inference problem. It is obvious that a particle filter would fail to recover the posterior distribution which is sharply centered about the observed value. We found that increasing the ensemble size to $M = 10000$ allows a particle filter to recover the posterior distribution; but the effective ensemble size still drops dramatically. The approximation provided by the EnKF is also displayed. The EnKF fails to recover the posterior due to its inherent linear regression ansatz which is inappropriate for this strongly nonlinear inference problem even in the limit ensemble size $M\to \infty$.

\begin{figure}[ht]
    \centering
    \includegraphics[width=0.48\textwidth]{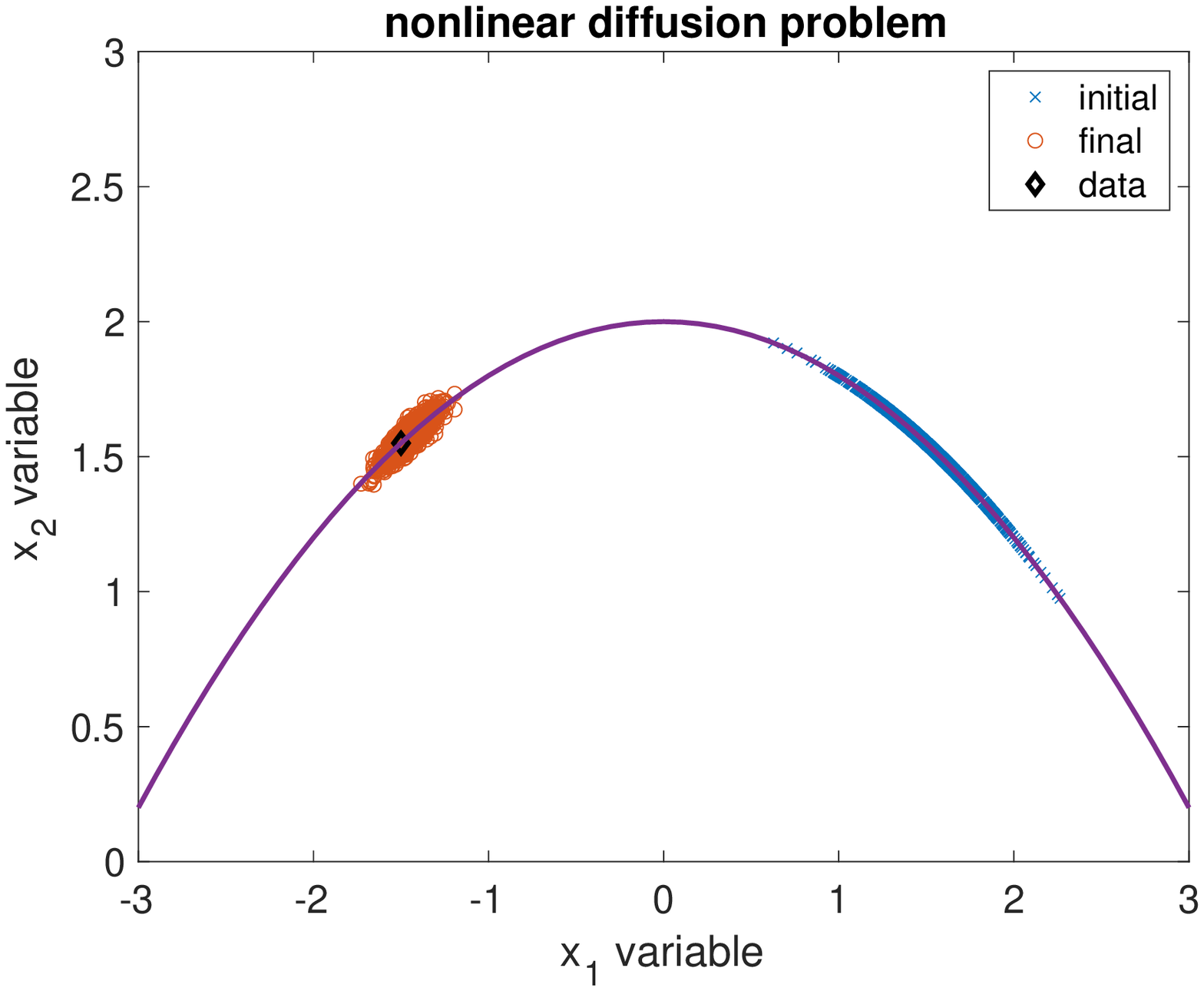}$\quad$
    \includegraphics[width=0.48\textwidth]{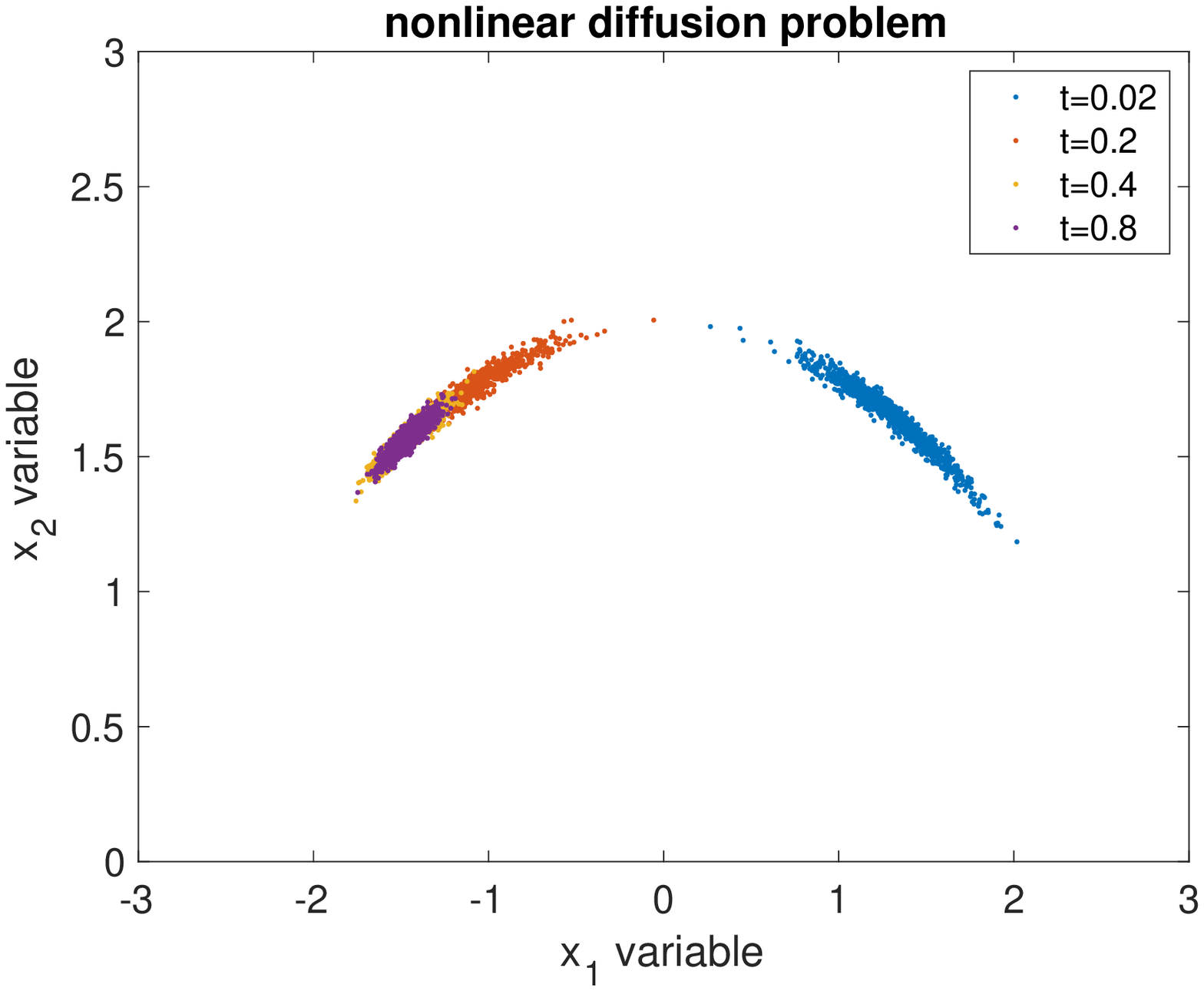}
    \caption{
        Left panel: Initial and final particle positions under the controlled evolution process. Right panel: Particle positions at intermediate times $t_k \in [0,1]$.
    }
    \label{fig:nonlinear2}
\end{figure}

In Figure \ref{fig:nonlinear2}, the results from the controlled mean field formulation are displayed. It can be concluded that the posterior distribution is will approximated despite the constant gain approximation made in order to formulate the control term $\hat g_t^{\rm KF}$ in (\ref{control_constant_gain}).

%
\subsection{Lorenz-63 example}
%

For a more challenging example we use the standard Lorenz-63 model \cite{lorenz63}
\begin{equation} \label{eq:lorenz}
    \frac{{\rm d}}{{\rm d}t} X_t = f(X_t) 
\end{equation}
where $X_t:\Omega \to  \mathbb{R}^3$ and
\begin{equation*}
f(x,y,z) = \left(\begin{array}{c}
           a(y-x) \\
           x(b-z)-y \\
           xy-cz
         \end{array} \right)
\end{equation*}
with parameters $a=10$, $b=28$ and $c=8/3$.

In order to obtain a reference solution $X_t^\dagger$ for $t\ge 0$, the ODE (\ref{eq:lorenz}) is solved numerically with
step-size $\Delta t = 0.005$ and initial condition
\[
X_0^\dagger = \left( \begin{array}{c} -0.587276 \\ -0.563678 \\16.8708 \end{array} \right).
\]
Scalar-valued observations are generated every $\Delta t_{\rm obs}>0$ units of time using the forward model
\[
y_{n\Delta t_{\rm obs}} = H X_{n\Delta t_{\rm obs}} + \nu_n, \qquad n = 1,\ldots N,
\]
with measurement measurement errors $\nu_n \sim {\rm N}(0,1)$ and forward map $H = (1 \,0 \,0) \in \mathbb{R}^{1\times 3}$. 
We use $\Delta t_{\rm obs} \in \{0.05,0.1,0.12\}$ in our experiments and perform $N = 20,000$ assimilation cycles.

The initial ensemble $\{X_0^{(i)}\}_{i=1}^M$ is drawn from the Gaussian distribution ${\rm N}(X_0^\dagger, 0.01 I)$.  We employ 
multiplicative ensemble inflation which amounts to replacing the Lorenz-63 dynamics by
\begin{equation*}
    \frac{{\rm d}}{{\rm d}t} X_t^{(i)} = f(X_t^{(i)})  + \sigma_k (X_t^{(i)} - \hat \mu_t ), \qquad i = 1,\ldots,M,
\end{equation*}
with inflation factors
\[
\sigma_k = 0.025 k, \qquad k = 0,\ldots,9.
\]
Here $\hat \mu_t$ denotes the empirical mean of the ensemble $\{X_t^{(i)}\}_{i=1}^M$. These equations are combined with
the augmented evolution equations (\ref{pure_drift}) and solved numerically with step-size $\Delta t = 0.005$ and ensemble
sizes $M \in \{5,10,15\}$. 

We report the resulting root mean square errors 
\[
\mbox{RMSE} = \sqrt{ \frac{1}{3N} \sum_{n=1}^N \|\hat \mu_{n\Delta t_{\rm obs}} - X_{n\Delta t_{\rm obs}}^\dagger \|^2 },
\]
which are computed for each ensemble size $M$, observation interval $\Delta t_{\rm obs}$ and inflation factor $\sigma_k$. 
The results are displayed in Table \ref{tab1} where the smallest RMSE over the range of
inflation factors $\{\sigma_k\}_{k=0}^9$ is stated  for each $M$ and $\Delta t_{\rm obs}$. We also state the corresponding RMSEs from 
a standard ensemble square root filter implementation \cite{asch2016data,Evensenetal2022}. We find that the proposed homotopy approach outperforms the ensemble square root filter in terms of RMSE in all settings considered. The improvements increase for increasing observation intervals $\Delta t_{\rm obs}$. The homotopy approach also appears less sensitive to the ensemble size $M$. 

\begin{table}
\begin{center}
\begin{tabular}{|c|c|c|c|}
\hline
$M$/$\Delta t_{\rm obs}$ & 5 & 10 & 15 \\
\hline 
0.05 & 0.5712/0.5457 & 0.5620/0.5475 & 0.5659/0.5496\\
0.10 & 0.8466/0.7735 & 0.8171/0.7627 & 0.8229/0.7707\\
0.12 & 0.9606/0.8645 & 0.9515/0.8621 & 0.9375/0.8615 \\
\hline
\end{tabular}
\end{center}
\caption{Displayed are the RMSE for both a standard ensemble square root filter implementation (top) and our homotopy approach (bottom) in terms of ensemble sizes $M \in \{5,10,15\}$ and observation intervals $\Delta t_{\rm obs} \in \{0.05,0.1,0.12\}$. The homotopy based data assimilation method leads to significantly reduced RMSEs in all settings.}
\label{tab1}
\end{table}

%
\section{Conclusions}
\label{sec:Conclusions}
%

Devising alternative proposal densities has a long history in the context of sequential data assimilation and filtering. Here we have explored a computationally tracktable approach which combines the concept of Schr\"odinger bridges with a rather straightforward homotopy approach. A further key ingredient is the approximate solution of the arising PDEs in terms of constant gain approximations, which are also widely used within the EnKF community. Numerical examples indicate that the approach is viable and can overcome limitations of both standard sequential Monte Carlo as well as standard EnKF methods. This has been demonstrated for single assimilation steps as well as long-time data assimilation using the chaotic Lorenz-63 model with only the first component observed infrequently. It remains to be seen how the proposed methods behave for high dimensional stochastic processes.

\medskip

\noindent
{\bf Acknowledgment.} This research has been funded by the Deutsche Forschungsgemeinschaft (DFG)- Project-ID 318763901 - SFB1294. We thank Nikolai Zaki for earlier work on the topic of this paper. 

\bibliographystyle{abbrvnat}

%
\bibliography{references}
%



%
\appendix
\section{Derivation of Control Term Equation}
%

\label{sec:forcing_derivation}
For an evolution equation
\begin{equation}
    \label{in_eq}
    dX_t = f(X_t)dt + \sqrt{2\sigma}dW_t
\end{equation}
we get a Fokker-Planck equation
\begin{equation}
    \label{in_fp}
    \partial_t \pi_t = -\nabla  \cdot \left( \pi_t \left( f - \sigma \nabla \log \pi_t \right) \right).
\end{equation}
Now we modify \eqref{in_eq} by an additional drift term, i.e.
\begin{equation}
    \nonumber
    dX^{\rm h}_t = f(X^{\rm h}_t)dt + \tilde g_t(X^{\rm h}_t)dt + \sqrt{2\sigma}dW_t
\end{equation}
with $\tilde g_t:\mathbb{R}^{d_x} \to \mathbb{R}^{d_x}$.
In that case, we would get a Fokker-Planck equation for the new equation:
\begin{equation}
    \label{fp_with_g}
    \partial_t \pi^{\rm h}_t = -\nabla  \cdot \left( \pi^{\rm h}_t \left( f + \tilde g - \sigma \nabla \log \pi^{\rm h}_t \right) \right)
    = -\nabla  \cdot \left( \pi^{\rm h}_t \left( f - \sigma \nabla \log \pi^{\rm h}_t \right) \right) - \nabla \cdot (\pi^{\rm h}_t \tilde g_t)
\end{equation}
We can find $\tilde g_t$ in terms of known quantities as follows: we begin by taking the derivative of $\pi^{\rm h}_t$ with respect to time:
\begin{equation}
    \nonumber
    \partial_t \pi^{\rm h}_t = - \pi^{\rm h}_t \left( \frac{\dot{Z}_t}{Z_t} + \frac{L}{T} \right)
    + \frac{1}{Z_t} e^{- \frac{t}{T} L } \partial_t \pi_t
\end{equation}
Next we substitute \eqref{in_fp} for $\partial_t \pi_t$ and use $\pi_t = Z_t e^{\frac{t}{T}L} \pi^{\rm h}_t$
\begin{align*}
    \partial_t \pi^{\rm h}_t =& - \pi^{\rm h}_t \left( \frac{\dot{Z}_t}{Z_t} + \frac{L}{T} \right)
    - Z_t^{-1} e^{- \frac{t}{T} L } \nabla  \cdot \left( \pi_t \left( f - \sigma \nabla \log \pi_t \right) \right)  \\
    =& - \pi^{\rm h}_t \left( \frac{\dot{Z}_t}{Z_t} + \frac{L}{T} \right)
       - Z_t^{-1}e^{- \frac{t}{T} L } \nabla  \cdot
           \left( Z_t e^{\frac{t}{T}L}\pi^{\rm h}_t
                \left( f - \sigma \nabla \log Z_t e^{\frac{t}{T}L}\pi^{\rm h}_t
            \right)
        \right)  \\
        =& - \pi^{\rm h}_t \left( \frac{\dot{Z}_t}{Z_t} + \frac{L}{T} \right)
       - e^{- \frac{t}{T} L } \nabla  \cdot
           \left( e^{\frac{t}{T}L}\pi^{\rm h}_t
                \left( f - \sigma \frac{t}{T}\nabla L - \sigma \nabla \log \pi^{\rm h}_t
            \right)
        \right)  \\
    =& - \pi^h_t \left( \frac{\dot{Z}_t}{Z_t} + \frac{L}{T} \right)
        - \nabla \cdot \left(\pi^{\rm h}_t \left(f- \sigma \frac{t}{T} \nabla L - \sigma \nabla \log \pi_t^{\rm h} \right) \right) \\
        & \qquad 
        - \,\,\frac{t}{T} \pi^{\rm h}_t \nabla L \cdot \left( f
        -\sigma \frac{t}{T}\nabla L- \sigma \nabla \log \pi_t^{\rm h}\right) .
\end{align*}
Comparing with \eqref{fp_with_g} it follows that we require
\begin{equation}
    \nonumber
    \nabla \cdot (\pi^{\rm h}_t \tilde g_t) =
    \frac{1}{T} \pi^{\rm h}_t \left( L + t \nabla L \cdot
        \left( f - \sigma \frac{t}{T} \nabla L - \sigma \nabla \log \pi^{\rm h}_t
    \right)
    \right)
    + \pi^{\rm h}_t \frac{\dot{Z}_t}{Z_t} - \nabla \cdot \left( \pi^{\rm h}_t \sigma \frac{t}{T} \nabla L \right)
\end{equation}
For the $\dot{Z}_t$ term we have
\begin{align*}
        \frac{\dot{Z}_t}{Z_t}
        =& \frac{1}{Z_t} \int \partial_t e^{-\frac{t}{T}L(x)} \pi_t(x) dx
        = \frac{1}{Z_t} \int -\frac{L}{T} e^{-\frac{t}{T}L(x)} \pi_t(x)
            + e^{-\frac{t}{T}L(x)} \partial_t \pi_t(x) dx \\
        =& -\frac{1}{T} \mathbb{E} L
        - \frac{1}{Z_t} \int e^{-\frac{t}{T}L(x)} \nabla \cdot \left( \pi_t \left( f -\sigma \nabla \log \pi_t
        \right) \right) dx \\
        =& -\frac{1}{T} \mathbb{E} L
        - \frac{1}{Z_t} \int \frac{t}{T} e^{-\frac{t}{T}L(x)} \pi_t \nabla L \cdot \left( f -\sigma \nabla \log \pi_t \right) dx \\
        =& -\frac{1}{T} \mathbb{E} L
        - \frac{t}{T} \mathbb{E} \nabla L \cdot \left( f -\sigma \nabla \log \pi_t \right) \\
        =& -\frac{1}{T} \mathbb{E} L
        - \frac{t}{T} \mathbb{E} \nabla L \cdot \left( f -\sigma \nabla \log \pi^{\rm h}_t - \sigma \frac{t}{T} \nabla L \right) \\
\end{align*}
where the third equality follows from integration by parts and the expected value is with taken with respect to $\pi^{\rm h}_t$.
We finally note that
\begin{equation} \nonumber
    \tilde g_t(x) = g_t(x) - \sigma \frac{t}{T}\nabla L(x).
\end{equation}


%
\section{Ensemble Kalman filter approximations}
\label{sec:linear_moment_matching}
%

We provide details on the derivation of the EnKF-like approximation (\ref{control_constant_gain}) to the controlled mean field equation (\ref{control_sde1}) and the various simplifications that arise from assuming a linear forward map.

We first recall that a continuous time formulation of the EnKF for a generic likelihood function
$$
L(x) = \frac{1}{2} (\hat h(x)-y_T)^\top \hat R^{-1}(\hat h(x)-y_T)
$$
is provided by
$$
\frac{\rm d}{{\rm d}t} X_t = - \Sigma_t^{x\hat h} \hat R^{-1} \left(
\frac{1}{2}\left(\hat h(X_t) + \pi_t[h]\right) - y_T\right).
$$
Here $\pi_t$ denotes the law of $X_t$, $\pi_t[g]$ the expectation value of a function $g$ under $\pi_t$, and
$$
\Sigma_t^{x\hat h} = \pi_t\left[ (x-\pi_t[x]) (\hat h - \pi_t[h])^\top\right]
$$
is the covariance matrix between the state $x$ and the forward map $\hat h$.

Formal application of this approach to the two contributions to the likelihood function (\ref{mod_likelihood}) leads to (\ref{control_constant_gain}). More precisely, the first term leads to $\hat h = h$ and 
$$
\hat R = \frac{t+\Delta t}{2\Delta t T}R
$$
while the second term results in $\hat h = \tilde h$ and
$$
\hat R = -\frac{t}{\Delta t T}R.
$$

The EnKF makes use of statistical linearization
$$
\Sigma_t^{xx} \pi_t[\nabla h] = \Sigma_t^{xh}
$$
which holds provided $\pi_t$ is Gaussian or if $h$ is linear; a result known as Stein's identity. The identity can also be used to approximate derivatives in a (weakly) non-Gaussian setting giving rise to
\begin{equation} \label{Stein}
\nabla h(x) \approx (\Sigma_t^{xx})^{-1} \Sigma_t^{xh}.
\end{equation}

We also recall the robust time-stepping method
\begin{equation} \label{time stepping}
X_{t_{n+1}}-X_{t_n} = -\Delta t 
\Sigma_{t_n}^{x\hat h}\left(\Delta t \Sigma_{t_n}^{\hat h\hat h}+
R\right)^{-1}\left(
\frac{1}{2}\left(\hat h(X_{t_n}) + \pi_{t_n}[h]\right) - y_T\right)
\end{equation}
which again can be adjusted appropriately to (\ref{control_constant_gain}).

We now assume a linear forward map, that is $h(x) = Hx$, and discuss the simplifications that result in the computation of (\ref{control_constant_gain}).
Note that
$$
\tilde h(x) = Hx - \Delta t Hf(x) + \Delta t \frac{\sigma t}{T}
H H^\top R^{-1} H x.
$$
Hence the covariance matrix $\Sigma_t^{x\tilde h}$ can be reformulated
to
$$
\Sigma_t^{x\tilde h} = \Sigma_t^{xx}H^\top - \Delta t
\Sigma_t^{xf}H^\top + \Delta t \frac{\sigma t}{T}
\Sigma_t^{xx}H^\top R^{-1} H H^\top
$$
and (\ref{control_constant_gain}) simplifies to
\begin{align*}
    \hat g_t^{\rm KF}(x) &=
    -\frac{1}{T} \Sigma_t^{xx}H^\top R^{-1}
    \left( \frac{1}{2} \left( Hx + H\mu_t^{\rm h}\right) - y_T \right)
    - \frac{t}{T} \Sigma_t^{xf}H^\top R^{-1} \left( \frac{1}{2} \left( Hx + H\mu_t^{\rm h}\right) - y_T \right) \\
    &\qquad +\,\,\frac{\sigma t^2}{T^2} \Sigma_t^{xx}H^\top R^{-1}
    H H^\top R^{-1} \left( \frac{1}{2} \left( Hx + H\mu_t^{\rm h}\right) - y_T \right)
    - \frac{t}{2T}\Sigma_t^{xx}H^\top R^{-1}H\left( 
    f(x) + \pi_t^{\rm h}[f]\right) \\
    & \qquad + \,\,\frac{\sigma t^2}{T^2}\Sigma_t^{xx} H^\top R^{-1}
    H H^\top R^{-1} \left( \frac{1}{2} \left(Hx+H\mu_t^{\rm h}\right) - y_T \right)
    + \mathcal{O}(\Delta t)\\
    &= -\left\{ \frac{1}{T}\Sigma_t^{xx} + \frac{t}{T}\Sigma_t^{xf}
    - \frac{2\sigma t^2}{T^2} \Sigma_t^{xx}H^\top R^{-1} H
    \right\}
    H^\top R^{-1} \left( \frac{1}{2} \left( Hx + H\mu_t^{\rm h}\right) - y_T \right)
    \\
    & \qquad -\,\, \frac{t}{2T}\Sigma_t^{xx}H^\top R^{-1}H\left( 
    f(x) + \pi_t^{\rm h}[f]\right) + \mathcal{O}(\Delta t) .
\end{align*}
Upon dropping terms of order $\mathcal{O}(\Delta t)$ and using $\Sigma_t^{xx}= \Sigma_t^{\rm h}$, we obtain (\ref{pure_diffusion}) for $f=0$ and (\ref{pure_drift}) for $\sigma = 0$ as special cases. The mean field equations (\ref{linear_homotopy}) also follow easily from $f(x) = Fx + b$.



\end{document}